\documentclass[12pt]{siamltex} 
\usepackage{epsfig, fullpage, cite, algorithm,algorithmic}

\def\frechet{Fr\'echet\ }

\def\eqnok#1{\/{\rm (\ref{eq:#1})}}
\def\begeq{\begin{equation}}
\def\endeq{\end{equation}}

\def\dnu{\, d \nu}
\def\dmu{\, d \mu}
\def\dlam{\, d \lambda}

\def\dint{{\displaystyle \int}}

\def\dfrac#1#2{\displaystyle \frac{#1}{#2}}

\def\calk{{{\cal K}}}

\newcommand{\lmin}{\lambda_{\mathrm{min}}}
\newcommand{\bx}{\bar x}
\newcommand{\by}{\bar y}
\newcommand{\gap}{\mathrm{gap}}

\title{
Condition Estimates for
Pseudo-Arclength Continuation
\thanks{Version of \today.
}
}
\author{K. I. Dickson%
\thanks{
North Carolina State University,
Center for Research in Scientific Computation and
Department of Mathematics,
Box 8205, Raleigh, N. C. 27695-8205, USA
({kidickso@unity.ncsu.edu, Tim\_Kelley@ncsu.edu, ipsen@math.ncsu.edu},
The work of these authors has been partially supported by National
Science Foundation Grants DMS-0404537 and DMS-0209695, and Army Research
Office Grants DAAD19-01-1-0592, W911NF-04-1-0276, and W911NF-05-1-0171.
}
\and\
C. T. Kelley\footnotemark[1]%
\and\
I. C. F. Ipsen\footnotemark[1]
\and\
I. G. Kevrekidis%
\thanks{
Department of Chemical Engineering,Princeton University,
Princeton, NJ 08544 ({yannis@princeton.edu}). Work supported in
part by AFOSR and an NSF/ITR grant.}
}

\begin{document}\maketitle

\begin{abstract}

We bound the condition number of the Jacobian in pseudo arclength
continuation problems, and we quantify the effect of this condition
number on the linear system solution in a Newton GMRES solve.

In pseudo arclength continuation one repeatedly solves systems of
nonlinear equations $F(u(s),\lambda(s))=0$ for a real-valued function
$u$ and a real parameter $\lambda$, given different values of the
arclength $s$. It is known that the Jacobian $F_x$ of $F$ with respect
to $x=(u,\lambda)$ is nonsingular, if the path contains only regular
points and simple fold singularities.  We introduce a new characterization
of simple folds in terms of the singular value decomposition, and we
use it to derive a new bound for the norm of $F_x^{-1}$. We also show
that the convergence rate of GMRES in a Newton step for
$F(u(s),\lambda(s))=0$ is essentially the same as that of the original
problem $G(u,\lambda)=0$. In particular we prove that the bounds on
the degrees of the minimal polynomials of the Jacobians $F_x$ and
$G_u$ differ by at most 2.  We illustrate the effectiveness of our
bounds with an example from radiative transfer theory.

\end{abstract}

\begin{keywords}
Pseudo-Arclength Continuation, singularity, GMRES, singular vectors, 
eigenvalues, rank-one update
\end{keywords}

%
%
\begin{AMS}
65H10, 
65H17, 
65H20, 
65F10, 
65F15  
\end{AMS}

\section{Introduction} 
\label{sec:realintro} Numerical continuation is the process of solving systems
of nonlinear equations $G(u,\lambda)=0$ for various values
of a real parameter $\lambda$. Here $u: R^N\rightarrow R$ 
is a real-valued function and $G: R^{N+1}\rightarrow R^N$.
An obvious approach for implementing numerical continuation, called
\textit{parameter continuation} \cite{herb, govaerts, kuznetsov},
traces out a solution path by repeatedly
incrementing $\lambda$ until the desired value of $\lambda$
is reached. In each such iteration, the current solution $u$
is used as an initial iterate for the next value of $\lambda$.
Although parameter continuation is simple and intuitive, it fails at
points $(u,\lambda)$ where the Jacobian $G_u$ is singular.
In this paper we consider singularities which are simple folds.

The standard way to remedy the failure of parameter continuation at
simple folds is to reparameterize the problem by introducing the
arclength parameter, $s$, so that both $u$ and $\lambda$ depend on
$s$.  This idea, known as \textit{pseudo-arclength continuation}
\cite{herb, govaerts, kuznetsov}, implements parameter continuation on $F(u(s),\lambda(s))=0$ with
$s$ as the parameter instead of solving $G(u,\lambda)=0$ with
$\lambda$ as the parameter.  Thus pseudo-arclength continuation requires
that the Jacobian $F_x$ of $F$ be nonsingular.  It is known that $F_x$
is nonsingular at simple folds and points where $G_u$ is nonsingular
\cite{herb}. 

Our first goal (\S~\ref{sec:nonsing}) is to quantify this nonsingularity. 
To this end we provide
a new characterization of simple folds in terms of the singular value
decomposition (SVD) of $G_u$. From the SVD, we derive a new bound for
$\|F_x^{-1}\|_2$.  This bound can be used to limit the arclength
step in Newton's method. As a byproduct we obtain a refinement
of Weyl's monotonicity theorem \cite{Par80} for the smallest eigenvalue of a 
symmetric positive semi-definite matrix (\S\ref{sec:supporttheorem}).

We also examine (\S\ref{sec:cluster})
how the conditioning of $F_x$ affects the convergence of the 
inner GMRES \cite{gmres} iteration in a Newton-GMRES  solver
\cite{ctk:roots,ctk:newton,brown/saad94,brown/saad90}.
We show that the eigenvalue clustering of the Jacobian $F_x$ in
the reformulated problem is not much different from that of the Jacobian
$G_u$ in the original problem. 
To be precise,  the upper bound on the degree of the minimal polynomial
of $F_x$ differs from that of $G_u$ by at most two.
This implies \cite{ctk:jordan,ctk:xue3} that the
convergence of GMRES as a linear solver in a  Newton step does not
slow down when parameter continuation is replaced by 
pseudo-arclength continuation.

At last (\S\ref{sec:heq}) we illustrate our findings with a numerical
example from radiative transfer theory.  These numerical results,
combined with our bounds, support the use of pseudo-arclength
continuation in solution paths that contain simple folds.

\section{Background}
\label{sec:intro}
We briefly review theory and algorithms for solving numerical
continuation problems $G(u,\lambda) = 0$, where $\lambda \in R$, 
$u: R^N \rightarrow R$ and $G: R^{N+1} \rightarrow R^N$.  We discuss parameter
continuation \S\ref{sec:paramc} and pseudo-arclength continuation
in \S\ref{sec:arcs}. We use the abbreviations
\[G_u \equiv \frac{\partial G}{\partial u}, \qquad
 G_{\lambda} \equiv \frac{\partial G}{\partial \lambda}.\]

\subsection{Simple parameter continuation}\label{sec:paramc}
Parameter continuation \cite{herb, govaerts, kuznetsov}
is the simplest method for solving $G(u,\lambda)=0$. The idea is to start at a point
$\lambda = \lambda_{init}$ and solve $G(u,\lambda)$ for $u$, say by
Newton's method. Use the solution $u_0$ as the initial iterate 
to solve the next problem  $G(u, \lambda + \dlam) = 0$.
Algorithm~{\tt paramc} below is a simple 
implementation of parameter continuation
from $\lambda_{init}$ to $\lambda_{end} = \lambda_{init} + n \dlam$ 
where $n$ denotes the maximum number of continuation iterations.

\bigskip

\begin{algorithm}
{$\mbox{\bf paramc}(u, G, \lambda_{init}, \lambda_{end}, \dlam)$}
\begin{algorithmic}
\STATE Set $\lambda = \lambda_{init}$, $u_0 = u$
\WHILE{$\lambda \le \lambda_{end}$}
\STATE Solve $G(u, \lambda) = 0$ with $u_0$ as the initial iterate to
obtain $u_1$
\STATE $u_0 = u_1$
\STATE $\lambda =\lambda + \dlam$
\ENDWHILE
\end{algorithmic}
\end{algorithm}

\bigskip

While parameter continuation appears to be a reasonable method for
solving $G(u,\lambda)=0$, it fails at points that violate the
assumptions of the implicit function theorem.  Such points of failure
are called singular points.

\medskip

\begin{definition}
\label{def:singpt}
A {\bf singular point} is a solution $(u_0, \lambda_0)$ to $G(u,\lambda)=0$ for which $G_u(u_0, \lambda_0)$ is singular.
\end{definition}

\medskip

In order to understand why parameter continuation fails at singular
points, we recall the implicit function theorem \cite{herb,ortega}. 
The norm $\| \cdot \|$
denotes the Euclidean norm, and $C^k (\Omega)$ denotes the space of
$k$ times continuously differentiable functions from an open subset
$\Omega \subset R^{N+1}$ to $R^N$.

\bigskip

\begin{theorem} 
\label{th:implicit}
[{\bf Implicit Function Theorem:}]
Let $\Omega$ be an open subset of $R^{N+1}$ and let $G \in C^k (\Omega)$
for some integer $k > 0$. Let $G_u$ and $G_\lambda$ be
Lipschitz continuous in $\bar \Omega$, the closure of $\Omega$.  If
\begin{itemize}
\item $(u_0, \lambda_0) \in \Omega$,
\item $G(u_0, \lambda_0) = 0$,
\item $G_u(u_0, \lambda_0)$ is nonsingular,
\end{itemize}
then there are $\rho>0$ and $\epsilon>0$ such that there is a unique solution
\[
v \in {\cal B}_{\rho}(u_0) \equiv \{ u \, | \, \| u - u_0 \| < \rho \}
\]
of $G(u, \lambda) = 0$ for all 
$\lambda \in (\lambda_0 - \epsilon , \lambda_0 + \epsilon)$
and $u \in {\cal B}_{\rho}(u_0)$.  Furthermore, $v$ is a $k$ times continuously differentiable function of $\lambda$.  
\end{theorem}

\medskip

If the assumptions of the implicit function theorem are satisfied then
Newton's method converges q-quadratically, as shown below.  We will use
the following definition of quadratic convergence for Newton's method.

\medskip

\begin{definition}
\label{def:qquad}
Let $\{ x_n \} \subset R^N$ be a sequence and let $x^* \in R^N$.  
We say that $x_n \rightarrow x^*$ {\bf q-quadratically} 
as $n\rightarrow \infty$, 
if $x_n \rightarrow x^*$ and if there is $K>0$ such that 
\[\Vert x_{n+1}-x^* \Vert \leq K\Vert x_n - x^* \Vert^2.\]   
\end{definition}

\medskip

The following corollary presents conditions under which
Newton's method applied to $G(u,\lambda)=0$ converges q-quadratically.

\medskip

\begin{corollary}
\label{cor:kant}
Let the assumptions of Theorem~\ref{th:implicit} hold. Then there is
$\delta>0$, which depends only on $\| G_u^{-1}(u_0,\lambda_0) \|$ and
the Lipschitz constants of $G_u$ and $u$, such that if $| \lambda -
\lambda_0 | < \delta$ then Newton's method with initial iterate $u_0$
converges q-quadratically to the solution $u^*$.
\end{corollary}

\smallskip

\begin{proof}  Define the Lipschitz constants 
$$\|u(\lambda)-u(\mu)\| \leq \gamma_u |\lambda-\mu|,\qquad 
\|G_u(u,\lambda)-G_u(v,\mu)\| \leq \gamma_G ( \| u-v| \| + |\lambda-\mu|).$$
According to \cite{ctk:roots} 
\[ \|u_0-u^*\| < \frac{1}{2\gamma_G\|G_u^{-1}(u_0,\lambda_0)\|} \]  so choosing 
\[ \delta < \frac{1}{2\gamma_u\gamma_G \|G_u^{-1}(u_0,\lambda_0)\|} \] completes the proof.
\end{proof}
 
\medskip 

Theorem~\ref{th:implicit} and Corollary~\ref{cor:kant}
suggest that points $(u,\lambda)$ for which $G_u$ is singular may
cause the loss of uniqueness in the solution to $G(u,\lambda)=0$ as well as
the failure of Newton's method, and therefore failure of
Algorithm~{\tt paramc}.  Thus we need a continuation method that
does not fail at singular points.

\bigskip

\subsection{Pseudo-arclength Continuation}
\label{sec:arcs}

Pseudo-arclength continuation \cite{herb,govaerts,kuznetsov}
avoids the  problems of Algorithm~{\tt paramc} at singular points
by using an arclength parameterization.
The curve in Figure~\ref{fig:heqdiag}, for instance, has a
singularity with respect to the parameter $\lambda$.  If we 
choose arclength $s$ as the parameter $\lambda$, and $x = (u^T,\lambda)^T$ 
in place of $u$, we can compute the curve with simple parameter continuation.
The curve in Figure~\ref{fig:heqdiag} has a simple fold, which is the
singularity of interest for this paper. Formally, a simple fold is
defined as follows.

\medskip

\begin{definition}
\label{def:simplefold}
A solution $(u_0, \lambda_0)$ of $G(u,\lambda) = 0$
is a {\bf simple fold} if
\begin{itemize}
\item $dim(Ker(G_u(u_0,\lambda_0)))=1$ and
\item $G_\lambda(u_0,\lambda_0) \not \in Range(G_u(u_0,\lambda_0))$.
\end{itemize}
\end{definition}

\medskip

To develop a pseudo-arclength continuation method, we assume that $x$ 
depends smoothly on $s$.  Then one can differentiate $G(u,\lambda) = 0$ 
with respect to $s$ and obtain
\begeq
\label{eq:dotxresid}
\frac{d G(u(s),\lambda(s))}{ds} = G_u {\dot u} + G_\lambda {\dot \lambda} = 0.
\endeq
Equivalently, one can differentiate $G(x) = 0$ and obtain $G_x {\dot x} = 0$.
Here, $\dot x$ denotes the derivative with respect to $s$.
Because the norm is the Euclidean norm and $s$ is arclength,
\begeq
\label{eq:arclength}
\| \dot x \|^2 = \| \dot u \|^2 + | \dot \lambda |^2 = 1.
\endeq
Since we introduced a new parameter $s$, we must add an equation to
$G(u,\lambda)=0$ so that the number of equations equals the number of
unknowns and we have a chance of obtaining a nonsingular Jacobian for
the reformulated problem.  Hence, we work with the extended
equations \begeq
\label{eq:fat}
F(x, s) = \left(
\begin{array}{c}
G(x) \\
{\cal N}(x, s)
\end{array}
\right) = \left( \begin{array}{c} 0 \\ 0 \end{array} \right).
\endeq
The normalization equation
${\cal N} = 0$ is an approximation of \eqnok{arclength} where
\begeq
\label{eq:norm1}
{\cal N}(x,s) = {\dot x}_0^T( x - x_0) - (s - s_0) = 0.
\endeq
This equation says that the new point on the path lies on the 
tangent vector through the current point $x_0$.

Given a known point $(x_0,s_0)$, the pseudo-arclength continuation method 
increments arclength by $ds=s - s_0$, and
solves \eqnok{fat} with the normalization \eqnok{norm1} by Newton's method
with initial iterate $x_0$.
Algorithm~{\tt psarc} is a simple implementation of pseudo-arclength
continuation.


\bigskip

\begin{algorithm}
{$\mbox{\bf psarc}(u, F,  s_{end}, \dlam)$}
\begin{algorithmic}
\STATE Set $s=0$,  $x_0 = (u_0^T, \lambda_0)^T$
\WHILE{$s \leq s_{end}$}
\STATE Approximate $\dot x$
\STATE Solve $F(x, s) = 0$ with $x_0$ as the initial iterate obtain $x_1$
\STATE $x_0 = x_1$
\STATE $s= s + ds$
\ENDWHILE
\end{algorithmic}
\end{algorithm}

\smallskip

Since pseudo-arclength continuation is just simple parameter 
continuation applied to $F$ with $s$ as the parameter, 
Corollary \ref{cor:kant} gives conditions for 
the convergence of Newton's method in pseudo-arclength continuation.

\begin{corollary}
Let the assumptions of Theorem~\ref{th:implicit} hold for $F$. Then there is
$\delta>0$, which depends only on $\| F_x^{-1}(x_0,s_0) \|$ and
the Lipschitz constants of $F_x$ and $x$, such that if 
$|s -s_0 | < \delta$ then Newton's method with initial iterate $x_0$
converges q-quadratically to the solution.
\end{corollary}

The proof of Corollary \ref{cor:kant} shows that the step
in arclength is bounded by
$$\delta < \frac{1}{2\gamma_x\gamma_F \|F_x^{-1}(x_0,s_0)\|},$$
where $\gamma_x$ and $\gamma_F$ are Lipschitz constants for $x$
and $F$, respectively. Therefore a bound on $\|F_x^{-1}\|$ 
is an important factor in bounding the arclength step.
In the next section we present the main result of this paper,
a new bound on $\|F_x^{-1}\|$.

\section{Nonsingularity of $F_x$}
\label{sec:nonsing}
For a solution $x_0=(u_0, \lambda_0)$ to $G(u,\lambda)=0$, we present an upper bound on $\|F_x^{-1}(x_0,s_0)\|$ in the case that
\begin{itemize}
\item $G_u(u_0,\lambda_0)$ is nonsingular or
\item $(u_0, \lambda_0)$ is a simple fold of $G(u,\lambda)=0$.
\end{itemize}
In order to derive the bound, we introduce a new characterization
of simple fold, which is based on the singular value decomposition
of $G_u$. We prove the bound in \S\ref{sec:proofmaintheorem}. 
In \S\ref{sec:supporttheorem} 
we refine Weyl's monotonicity theorem for the smallest eigenvalue of a 
symmetric positive semi-definite matrix, which we need for the proof.

Let
\[
G_u(u, \lambda) = U \Sigma V^T
\]
be a singular value decomposition (SVD) of $G_u(u,\lambda)$ where
\[
\Sigma = diag(\sigma_1, \sigma_2, \dots, \sigma_N), \qquad \sigma_1 \geq \sigma_2 \geq \dots \geq \sigma_N \quad and \qquad u_N \equiv U e_N
\]
where $e_N$ is the last column of the $N \times N$ identity matrix.
The trailing column $u_N$ of $U$ is a left singular vector associated 
with the smallest singular value $\sigma_N$.
Since the singular values are continuous functions of the elements 
in $G_u(u,\lambda)$, they are also continuous in $\lambda$. If
\[\sigma_{N-1} \ge {\bar \sigma} > 0\]
for all $(u,\lambda)$ then the nullity of $G_u(u,\lambda)$ is at most one.  
If in addition $\sigma_N = 0$ then $u_N$ 
spans the left nullspace of  $G_u(u, \lambda)$.  From the direct sum
\[Ker(G_u^T) \oplus Range(G_u) = R^N\]
we see that $G_{\lambda}(u_0,\lambda_0)$ is not in the $Range(G_u)$
if and only if $G_{\lambda}^T u_N\neq 0$.
Hence we have a new, equivalent definition of simple fold.

\begin{definition}[Simple Fold via SVD]
Let $(u_0, \lambda_0)$ be a solution of $G(u,\lambda) = 0$,
and let $u_N$ be a left singular vector of $G_u(u_0,\lambda_0)$
associated with $\sigma_N$.

Then $(u_0,\lambda_0)$ is a {\bf simple fold} if
\begin{itemize}
\item $dim(Ker(G_u(u_0,\lambda_0)))=1$ and
\item $u_N^TG_\lambda(u_0,\lambda_0)\neq 0$.
\end{itemize}
\end{definition}

Continuity of $G_\lambda^T u_N$ implies that there is 
$\alpha > 0$ such that for all $(u, \lambda)$
\[\max\left( \sigma_N^2, |u_N^T G_\lambda|^2\frac{\gap}{\gap + \xi^2}
\right) \ge \alpha > 0,\]
where
\[
\gap \equiv \sigma_{N-1}^2 - \sigma_{N}^2, \qquad\mbox{and}\quad
\xi \equiv | u_N^T G_\lambda | + \| (I - u_N u_N^T) G_\lambda \|.
\]

\newcommand{\smin}{\sigma_{\mathrm{min}}}

\begin{theorem}\label{maintheorem}
Let $\bar{\Omega}$ be the closure of an open subset $\Omega \in R^{N+1}$,
and let $G$ be continuously differentiable in $\bar{\Omega}$.
Let $x_0=(u_0,\lambda_0)$ in $\bar{\Omega}$ be a solution to $G(u_0,\lambda_0)=0$,
and ${\cal N}(x_0,s_0)=0$ with $\|{\dot x}_0\|=1$. 
Let $\tau\geq 0$ be such that 
$\|G_u{\dot u}_0+G_{\lambda}{\dot\lambda}_0\|\leq \tau$.

Assume that for all $(u,\lambda)$ in ${\bar\Omega}$ there exists 
$\alpha>0$ such that 
$$\sigma_{N-1}>0,\qquad 
\max\left\{\sigma_N^2,\,(u_N^TG_{\lambda})^2\,{\gap\over \gap+\xi^2}\right\}
\geq \alpha,$$
where 
$$\gap\equiv \sigma_{N-1}^2-\sigma_N^2,\qquad 
\xi \equiv |u_N^TG_{\lambda}|+\|(I-u_N u_N^{T})G_{\lambda}\|.$$

If $\tau<\alpha$, then for all $x=(u,\lambda)$ in $\bar{\Omega}$, 
the smallest singular value $\smin(F_x)$ of the Jacobian $F_x$ of 
$F(x,s)$  is bounded from below with
$$\smin(F_x)\geq \sqrt{1-
\tau\max\left\{{1\over \alpha},1\right\}}.$$
\end{theorem}

We postpone the proof of Theorem \ref{maintheorem} 
until \S\ref{sec:supporttheorem} in order to derive an auxiliary result first.

\subsection{Lower Bound for the Smallest Eigenvalue}
\label{sec:supporttheorem}

We derive a lower bound for
the smallest eigenvalue of the rank-one update $A+yy^T$,
where $A$ is a real symmetric positive semi-definite matrix of order $N$, 
and $y$ is a real $N\times 1$ vector.

Let $\beta_1\geq \ldots\geq  \beta_N$ be the eigenvalues of $A$.
Weyl's monotonicity theorem \cite[Theorem (10.3.1)]{Par80} implies 
bounds for the smallest eigenvalue of $A+yy^T$:
$$\beta_N\leq \lmin(A+yy^T)\leq \beta_{N-1}.$$
Intuitively one would expect that $\lmin(A+yy^T)$ is larger
if $y$ is close to an eigenvector of $\beta_N$.
We confirm this by deriving lower bounds for $\lmin(A+yy^T)$
that incorporate the angle between $y$ and the eigenspace of $\beta_N$.

\begin{theorem}\label{l_1}
Let $A$ be an $N\times N$ real symmetric positive semi-definite matrix,
$u_N$ an eigenvector of $A$ associated with $\beta_N$,  $\|u_N\|=1$, and
$y\neq 0$ a real $N\times 1$ vector. Set $y_N\equiv u_N^Ty$. 
Then
\begin{equation}\label{eq:mainbound}
\lmin(A+yy^T)\geq \max\left\{\beta_N,\,
y_N^2\>{\gap\over \gap+\xi^2}\right\}
\end{equation}
where $\gap\equiv \beta_{N-1}-\beta_N$ and 
$\xi\equiv |y_N|+\sqrt{\|y\|^2-|y_N|^2}$.
\end{theorem}

\begin{proof}
We first show that
\begin{equation}\label{eq:boundhelper}
\lmin(A+yy^T)\geq \min\{\beta_N+y_N^2{\gap\over\gap+\xi^2},\>
\beta_{N-1}{y_N^2\over\xi^2}\}
\end{equation}
is lower bound for $\lmin(A+yy^T)=\min_{\|x\|=1}{x^T(A+yy^T)x}$.

Let 
$$A=U\pmatrix{\beta_1& & \cr &\ddots& \cr && \beta_N}U^T$$
be an eigendecomposition of $A$, and
$x$ be any real vector with $\|x\|=1$. Partition
$$U^Tx=\pmatrix{\bx \cr x_N},\qquad U^Ty=\pmatrix{\by\cr y_N}$$
so that $\xi=|y_N|+\|\by\|$.
Then 
$$x^T(A+yy^T)x\geq \beta_{N-1}\|\bx\|^2+\beta_Nx_N^2+(y^Tx)^2.$$
If $\|\bx\|\geq |y_N|/\xi$ then 
$$x^T(A+yy^T)x\geq (\beta_{N-1}y_N^2)/\xi^2,$$
which proves the second part of the bound in \eqnok{boundhelper}.

If $\|\bx\|<|y_N|/\xi$ then $|y_N|-\|\bx\|\xi>0$,
and it makes sense to use $|x_N|\geq 1-\|\bx\|$ in
$$|y^Tx|=|y_Nx_N+\by^T\bx|\geq |y_Nx_N|-\|\bx\|\|\by\|\geq
|y_N|-\|\bx\|\xi.$$
Hence
$$x^T(A+yy^T)x\geq \beta_{N-1}\|\bx\|^2+\beta_Nx_N^2+(y^Tx)^2
\geq\beta_N+y_N^2+(\gap+\xi^2)\|\bx\|^2-2\xi \|\bx\||y_N|.$$
This is a function of $\|\bx\|$ which has a minimum at 
$\|\bx\|=|y_N|\xi/(\gap+\xi^2)$.
Hence
$$x^T(A+yy^T)x\geq \beta_N+y_N^2{\gap\over \gap+\xi^2},$$
which proves the first part of the bound in \eqnok{boundhelper}.

With the help of \eqnok{boundhelper} we now show the desired bound 
\eqnok{mainbound}.  
Weyl's theorem \cite[Theorem (10.3.1)]{Par80} implies 
$\lmin(A+yy^T)\geq \beta_N$, which proves the first part of the bound in \eqnok{mainbound}.
For the second part of the bound in \eqnok{mainbound}, we use the fact that the eigenvalues
of $A$ are non-negative, hence $\beta_{N-1}\geq \gap$ and
$${\beta_{N-1}\over \xi^2}\geq {\gap\over \gap +\xi^2}.$$
Substituting this into \eqnok{boundhelper} gives the second part
of the bound in \eqnok{mainbound}
\begin{eqnarray*}
\min(A+yy^T)&\geq& \min\{\beta_N+y_N^2{\gap\over \gap+\xi^2},\,
y_N^2 {\beta_{N-1}\over \xi^2}\}\\
&\geq&
\min\{\beta_N+y_N^2{\gap\over \gap+\xi^2},\,y_N^2{\gap\over \gap+\xi^2}\}
=y_N^2{\gap\over \gap+\xi^2}.
\end{eqnarray*}

\end{proof}

The quantity $\gap$ in Theorem \ref{l_1} is the absolute gap
between the smallest and next smallest eigenvalues.
The theorem shows that $\lmin(A+yy^T)$ is likely to be larger if 
$y$ has a substantial contribution in the eigenspace of $\beta_N$.
The bound in Theorem \ref{l_1} is tight when $y$ is a multiple of $u_N$. 
That is, if $|u_N^Ty|=\|y\|$ then
$\lmin(A+yy^T)= \min\{\beta_N+\|y\|^2,\beta_{N-1}\}$.

Now we are in a position
to complete the proof of Theorem~\ref{maintheorem}.

\subsection{Proof of Theorem \ref{maintheorem}}\label{sec:proofmaintheorem}
Define the residual $r\equiv G_u{\dot u}_0+G_{\lambda}{\dot\lambda}_0$ and
form
$$F_xF_x^T=\pmatrix{G_u & G_{\lambda} \cr {\dot u}_0^T & {\dot\lambda}_0}
\pmatrix{G_u^T &{\dot u}_0 \cr G_{\lambda}^T & {\dot\lambda}_0}=
\pmatrix{G_uG_u^T+G_{\lambda}G_{\lambda}^T & r \cr r^T &1}.$$
The eigenvalues of $F_xF_x^T$ are the squares of the singular values
of $F_x$.
Applying Theorem~\ref{l_1} to $G_uG_u^T+G_{\lambda}G_{\lambda}^T$
with $A=G_uG_u^T$, $y=G_{\lambda}$, $\beta_N=\sigma_N^2$,
$\beta_{N-1}=\sigma_{N-1}^2$ and $\gap=\sigma_{N-1}^2-\sigma_N^2$
shows $\lmin(G_uG_u^T+G_{\lambda}G_{\lambda}^T)\geq \alpha$.
Hence we can write
$$\pmatrix{G_uG_u^T+G_{\lambda}G_{\lambda}^T & 0 \cr 0&1}^{-1}
F_xF_x^T=I+E,$$
where $\|E\|\leq \tau\max\left\{{1\over \alpha},1\right\}$.
If $\tau<\min\{\alpha,1\}$ then $\|E\|<1$, $I+E$ is nonsingular, and
$${1\over \|(F_xF_x^T)^{-1}\|}\geq 1-
\tau\max\left\{{1\over\alpha},1\right\}.$$


\section{Newton-GMRES and Eigenvalue Clustering}
\label{sec:cluster}

This section discusses the performance of the inner GMRES
iteration in the context of continuation with a Newton-GMRES nonlinear
solver.
Theorem~\ref{maintheorem} gives conditions under which the 
Jacobian matrix $F_x$ of the reformulated problem is uniformly nonsingular.  
This implies GMRES is a practical candidate for making the linear 
solve in Newton's method when implementing pseudo-arclength continuation.  
While the results in the previous section address
conditioning, they do not directly translate into the performance
of iterative methods \cite{ctk:roots,trefbau,IpsM97}, 
especially in the non-normal
case.  However, we can go further to see that the eigenvalue clustering properties of the matrix $F_x$ do not stray far from those of $G_u$.  

Suppose the eigenvalues of $G_u$ are nicely clustered (in the sense of
\cite{ctk:jordan,ctk:xue3}). Even in the singular case, this would mean
that the zero eigenvalue of $G_u$ is an ``outlier''. We seek to show
that adding the row and column does not significantly increase the number
of outliers, and that we can then use the estimates in
\cite{ctk:jordan,ctk:xue3}.

One approach is to use the paradigm of \cite{ctk:yannis1}.  The idea is
that
\begeq
\label{eq:gdecomp}
G_u = I + K(u) + E
\endeq
where $K_u$ is a low-rank operator, say of rank $p$, 
and $E$ is small. We then want to
write $F_x$ in the same way, and then compare the number of outliers
by comparing the ranks of the $K$-terms.

Assume that $E$ is small enough
so that the eigenvalues of $I - K$ are ``outliers'' in the sense of
\cite{ctk:jordan}. Since the degree of the minimal polynomial
of $I - K$ is at most $p+1$, we have a bound 
for the sequence of residuals $\{ r_l \}$ of the GMRES iteration 
of the form

\begeq \label{eq:jbound}
\| r_{{\hat p}+k} \| \le C \| E \|^k \| r_0 \|
\endeq
where ${\hat p} \le p+1$ GMRES iterations are needed to kill the
contribution of the outlying eigenvalues. 

Theorem~\ref{th:yannis} states that
that the spectral properties of $F_x$ are similar to those of $G_u$. 

\begin{theorem}
\label{th:yannis}
Let the assumptions of Theorem~\ref{maintheorem} hold. Assume that
\eqnok{gdecomp} holds with $rank(K(u)) = p$. Then there is $\calk(u)$
having rank at most $p+2$ such that
\[
\| F_x - I - \calk(u) \| \le \| E \|.
\]
\end{theorem}

\begin{proof}
We write
\cite{ctk:yannis1}
\[
F_x = I_{(N+1) \times (N+1)} +
\left(
\begin{array}{cc}
K & G_\lambda\\
{\dot u}^T & {\dot \lambda}
\end{array}
\right)
+
\left(
\begin{array}{cc}
E & 0 \\
0 & 0
\end{array}
\right).
\]
The range of
\[
{\cal K} = \left(
\begin{array}{cc}
K & G_\lambda\\
{\dot u}^T & {\dot \lambda}
\end{array}
\right)
\]
is

$$
\left(
\begin{array}{c}
Range(K)\\ 0
\end{array}
\right) +
\mbox{span } \left\{
\left(
\begin{array}{c}
G_{\lambda}\\ 0
\end{array}
\right)
\right\}+
\mbox{span} \left\{
\left(
\begin{array}{c}
0\\ 1
\end{array}
\right)
\right\}
$$
and hence the rank of $\cal K$ is at most $p+2$.
\end{proof}

So, while the eigenvalues may change, we have not increased the degree of
the minimal polynomial of the main term ($K$ vs $\cal K$) beyond $p+3$.
Hence, the methods of \cite{ctk:jordan} can be applied to obtain
a bound like \eqnok{jbound} with ${\hat p} \le p+3$.

\section{Example: Chandrasekhar H-Equation}
\label{sec:heq}

We now present an example of a solution path containing a simple fold.  The equation of interest is called the Chandrasekhar $H$-equation \cite{chand,ctk:roots,twm68} from radiative transfer theory:

\begeq\label{eq:heq1}
H(\mu) = 1 + \frac{c}{2} H(\mu) \dint_0^1 H(\nu) \dfrac{\dnu \mu}{\mu + \nu}.
\endeq
The goal is to compute the $l^1$ norm of the solution to Equation~\eqnok{heq1} for various natural parameter values $c$.  That is, we compute
\[
\| H \|_1 = \dint_0^1 H(\nu,c) \dnu
\]
as a function of $c$. Integrating ~\eqnok{heq1} with respect to $\mu$ yields
\[
\| H \|_1 = 1 + \frac{c}{2}
\dint_0^1 \dint_0^1 \dfrac{H(\mu) H(\mu) \mu \dmu \dnu}{\mu + \nu}
= 1 + \frac{c}{4} \| H \|_1^2,
\]
and so
\begeq
\label{eq:hnorm}
\| H \|_1 = \dfrac{1 \pm \sqrt{1 - c}}{c/2}.
\endeq

Equation~\eqnok{hnorm} tells us two interesting things. First, there can
be no real solutions of the $H$-equation for $c > 1$, so there
must be a singularity at $c=1$, or else the implicit function theorem would
tell us that we could continue past $c=1$. Secondly, the $\pm$ gives us
a hint that there may be two solutions, at least for $0 < c < 1$ (and
there are!).

\begin{figure}[h!]
\caption{\label{fig:heqdiag} $\|H \|_1$ as a function of $c$}
\centerline{
\epsfxsize=3.5in
\epsfbox{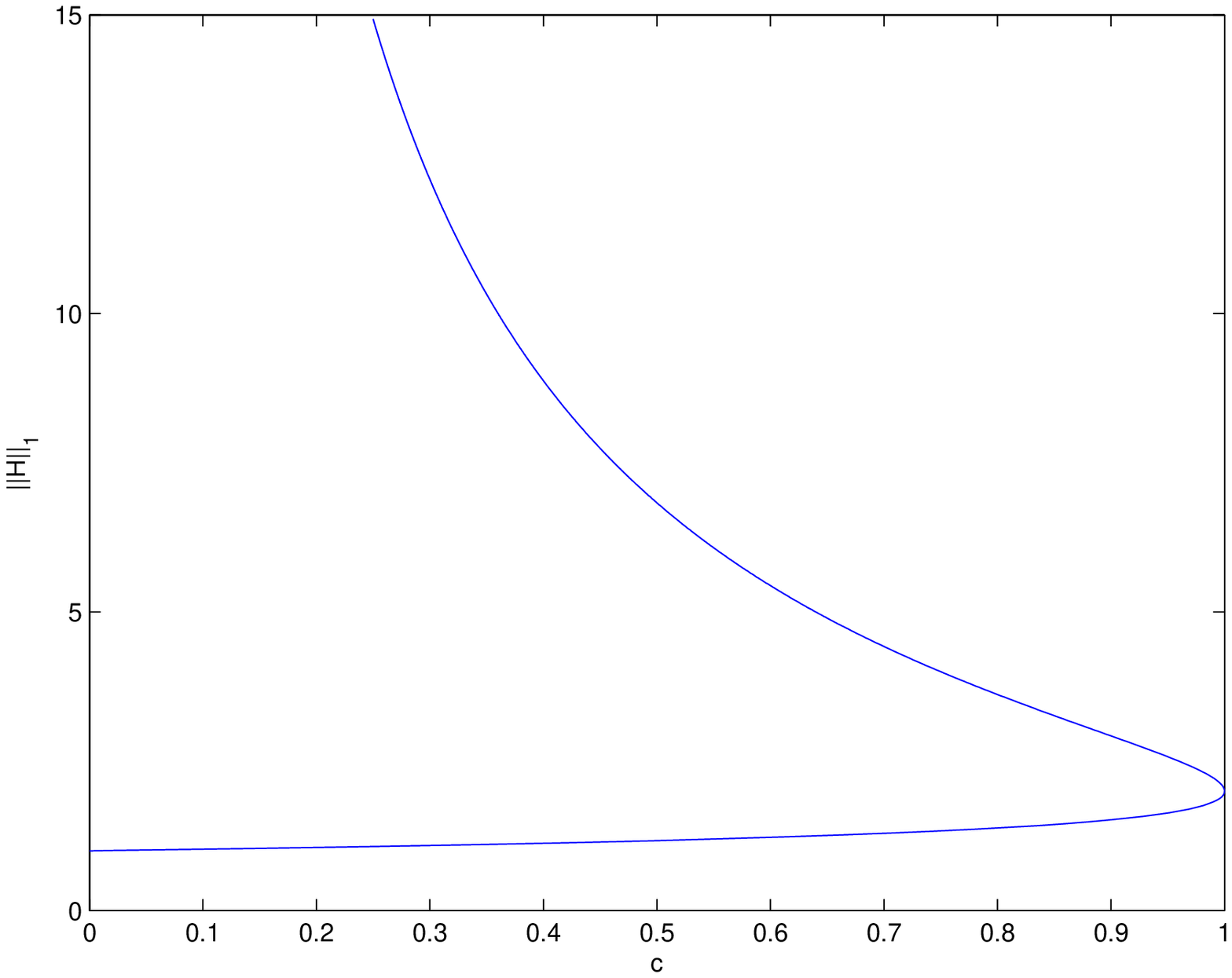}
\leavevmode
}
\end{figure}

Figure~\ref{fig:heqdiag} is a plot of $\| H \|_1$ against $c$.  Notice how the curve bends around when $c =1$, and how there are
two solutions for each $0 < c < 1$.  In fact, we are witnessing a simple fold at $c=1$.

\subsection{Simple Fold at $c=1$}  For the $H$-equation, it is possible to compute
the singularity analytically.  Write the $H$-equation as
\[
G(H,c)(\mu) = H(\mu) - \left( 1 - \frac{c}{2} \dint_0^1
\dfrac{\mu H(\nu) \dnu }{\mu + \nu}
\right)^{-1}.
\]
Taking the \frechet derivative of $F$ in the direction of $w$ yields
\[
G_H(H,c) w(\mu) = w(\mu)- \dfrac{
\frac{c}{2} \dint_0^1 \dfrac{\mu w(\nu) \dnu }{\mu + \nu} }%
{ \left( 1 - \frac{c}{2}
\dint_0^1 \dfrac{\mu H(\nu) \dnu }{\mu + \nu} \right)^2 }
= w(\mu) - \frac{c}{2} H(\mu)^2
\dint_0^1 \dfrac{\mu w(\nu) \dnu }{\mu + \nu}.
\]

Let $c=1$, then \eqnok{hnorm} implies that
\[
\dint_0^1 H(\mu) \dmu = 2
\]
and therefore
\[
\begin{array}{ll}
\frac{1}{2} \dint_0^1 \dfrac{\nu H(\nu) \dnu }{\mu + \nu}
& =
\frac{1}{2} \dint_0^1 H(\nu) \left( 1 - \dfrac{\mu}{\mu + \nu} \right)
\dnu 
= 1 - \dint_0^1 \dfrac{\mu H(\nu) \dnu }{\mu + \nu} \\
\\
& = H(\mu)^{-1}.
\end{array}
\]
Hence if $\phi(\mu) = \mu H(\mu)$,
\[
G_H(H,1) \phi = 0,
\]
and we have shown directly that $G_H$ is singular at $c=1$.

One can apply Perron-Frobenius theory
\cite{karlin,carl} to show that the null space of $G_H$ has
dimension one, and hence is spanned by $\phi$. The singularity at
$c=1$ is a simple fold because
\[
G_c(H,1)(\mu) = - H^2(\mu) \frac{1}{2} \dfrac{\mu H(\nu) \dnu }{\mu + \nu}
= H^2(\mu) (H^{-1}(\mu) - 1)   
\]
is not in the range of $G_H$. To see this note that 
\[
G_c(H,1)(\mu) = H^2(\mu) (H^{-1}(\mu) - 1) \le 0
\]
and vanishes only at $\mu = 0$. The null space of $G_{H}^{T}$ is the span of $H^{-1}$, which is 
strictly positive. Hence $G_c$ is not orthogonal to the null space of
$G_{H}^T$.

One can also show that $G_H$ is nonsingular for all $c \ne 1$ by 
an argument even more tedious than the one above \cite{ctk:hnewt}.

\subsection{Smallest Singular Values}  
\begin{figure}[h!]
\caption{\label{fig:smallsingvals} $\sigma_{\min}(F_{(H,c)})$ as a function of $c$}
\centerline{
\epsfxsize=3.5in
\epsfbox{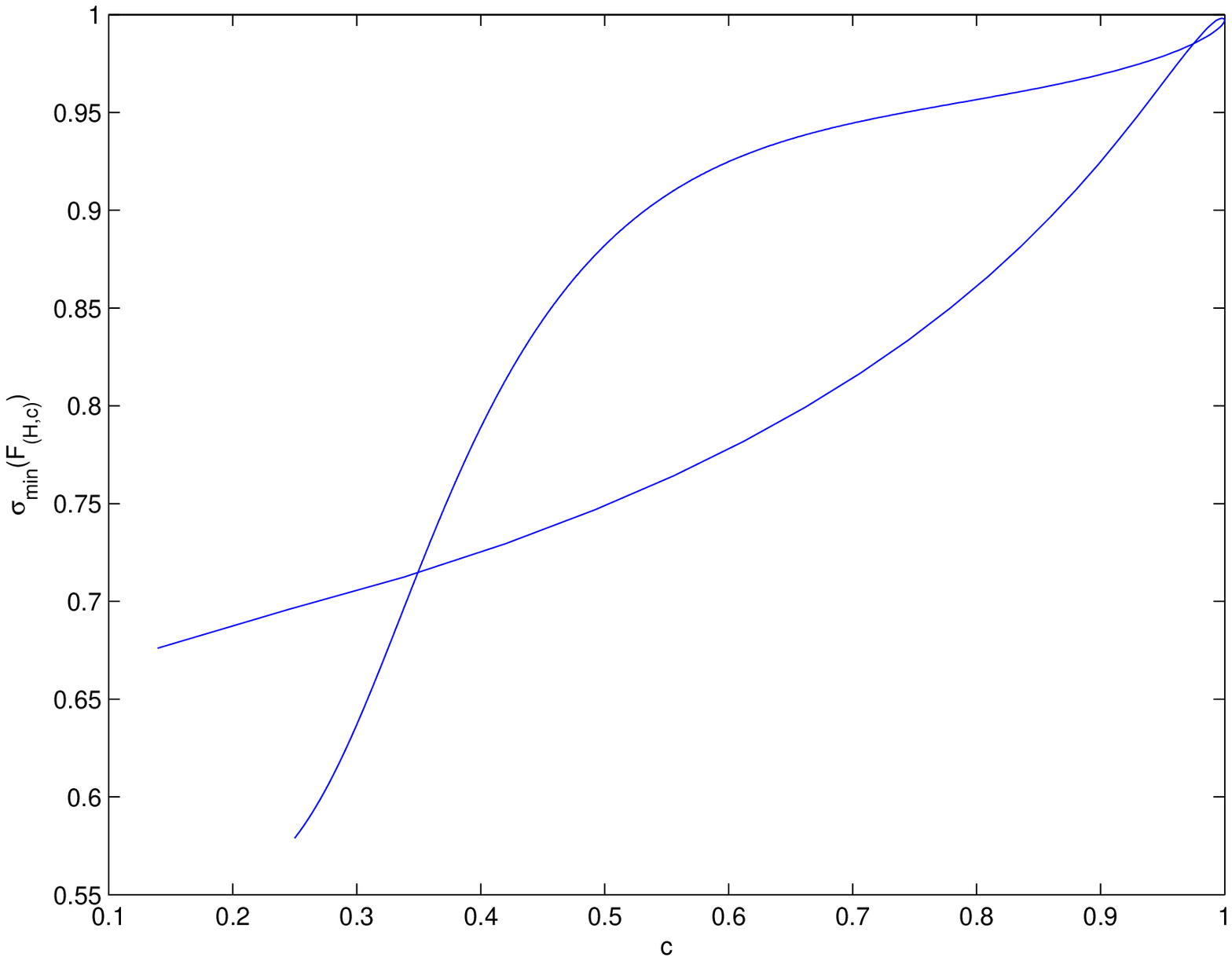}
\leavevmode
}
\end{figure}

As a demonstration of the result in \S~\ref{sec:nonsing}, we calculate the smallest singular value of the Jacobian matrix associated with the augmented system for the $H$-equation with each continuation iteration.  In the language of \S~\ref{sec:nonsing}, we find $\sigma_{\min}(F_{(H,c)})$ for various $c$ where $F_{(H,c)}$ denotes the Jacobian of $\left(\begin{array}{c} G(H,c) \\ {\cal N}(H,c,s) \end{array}\right)$.  Figure~\ref{fig:smallsingvals} shows that the smallest singular value of  $F_{(H,c)}$ for each $c$ stays away from zero keeping $F_{(H,c)}$ nonsingular, even at the simple fold ($c=1$).  The pseudo-arclength code used here uses a direct $LU$ factorization of the Jacobian for the linear solve in Newton's method.  The step in arclength is fixed at $ds=.5$, and we use a a secant predictor \cite{herb}.  The integral is discretized with the composite midpoint rule and 200 nodes.  The singular values are calculated using Matlab's {\tt svd} command.

\subsection{Computation with $H$-equation}

The consequences of the remarks in \S~\ref{sec:cluster} are that
for
a problem like the $H$-equation, which is a nonlinear compact fixed
point problem, the number of GMRES iterations per Newton step should
be bounded. One must take this expectation with a grain of salt because
as one moves along the path, the norm of the solution increases, and so
the number of outliers may increase slowly.  The observations we present
illustrate this.

We use a Newton-GMRES version of pseudo-arclength continuation
\cite{ctk:ferng1}, fixing the step in arclength to $ds = .02$, using
a secant predictor \cite{herb}, and beginning the continuation at $c=0$,
where the $H = 1$ is the solution.  The vector with coordinates all
equal to one is the solution of the discrete problem as well. We discretize
the integral with the composite midpoint rule using 400 nodes.

\begin{figure}
\caption{\label{fig:kpn} Krylov's per Newton}
\centerline{
\epsfxsize=3.5in
\epsfbox{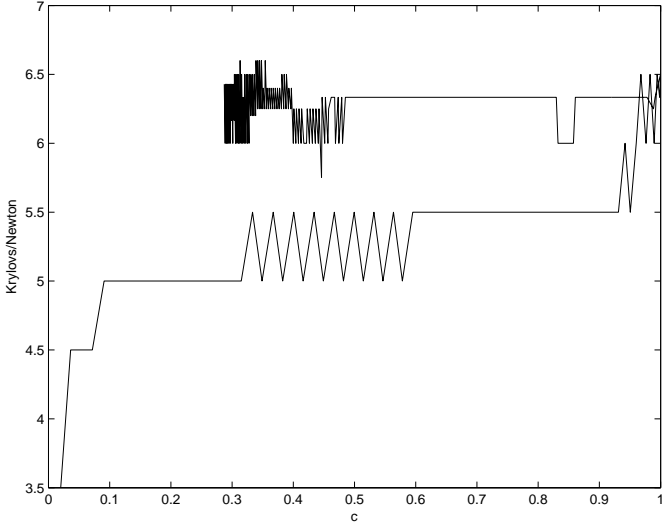}
\leavevmode
}
\end{figure}

In Figure~\ref{fig:kpn} we plot the average number of GMRES iterations per
Newton iteration as a function of $c$. As one moves further on the path, the
predictor becomes less effective, and the number of Newton iterations increase. Moreover, the norm of the solution also increases adding roughly
one to the number of Krylov's per Newton.

\clearpage


\vspace{-.25in}
%

\end{document}